\documentclass[12pt,titlepage]{article}
\usepackage{amssymb}
\input epsf
\newcommand{\proof}{{\noindent \bf Proof. }}
\newtheorem{thm}{Theorem}
\newtheorem{defi}[thm]{Definition}

\newtheorem{lem}[thm]{Lemma}
\newtheorem{cor}[thm]{Corollary}
\newtheorem{prop}[thm]{Proposition}

\def\2{\mathbb Z_2}

\def\qed{$\Box$}

\def\SG{{\rm SG}}

\title{On directed local chromatic number, shift graphs, and Borsuk-like
  graphs}  
\author{
      {\bf G\'abor Simonyi}\thanks{Research partially supported by the
Hungarian Foundation for Scientific Research Grant (OTKA) Nos.\ 
AT048826, T076088.}\\
 Alfr\'ed R\'enyi Institute of Mathematics\\
 Hungarian Academy of Sciences\\
 1364 Budapest, POB 127, Hungary\\
 {\tt simonyi@renyi.hu}
\and
   {\bf G\'abor Tardos}\thanks{Research partially supported by the
NSERC grant 611470 and the
Hungarian Foundation for Scientific Research Grant (OTKA) Nos.\ 
AT048826.}\\
 School of Computing Science \\
 Simon Fraser University \\
 Burnaby BC, Canada V5A 1S6\\
 and\\
 Alfr\'ed R\'enyi Institute of Mathematics\\
 Hungarian Academy of Sciences\\
 1364 Budapest, POB 127, Hungary\\
 {\tt tardos@cs.sfu.ca}}
\date{}

\begin{document}
\maketitle

\begin{abstract}
We investigate the local chromatic number of shift graphs and prove that it
is close to their chromatic number. This implies that the gap between the
directed local chromatic number of an oriented graph and the local chromatic
number of the underlying undirected graph can be arbitrarily large. We also
investigate the minimum possible directed local chromatic number of oriented
versions of ``topologically $t$-chromatic'' graphs. We show that this minimum
for large enough $t$-chromatic Schrijver graphs and $t$-chromatic generalized
Mycielski graphs of appropriate parameters is $\lceil t/4\rceil+1$. 
\end{abstract}

\section{Introduction}

The local chromatic number of a graph $G$, defined by Erd\H{o}s, F\"uredi,
Hajnal, Komj\'ath, R\"odl, and Seress \cite{EFHKRS} is a coloring parameter
that was further investigated recently in the papers \cite{KPS, ST,
  STV}. (See also \cite{BdW} for some related results.) Denoting the set of
neighbors of a vertex $v$ by $N(v)$, it is defined as follows.  

\begin{defi} \label{defi:lochr} {\rm(\cite{EFHKRS})}
The {\em local chromatic number} of a graph $G$ is
$$\psi(G):=\min_c \max_{v\in V(G)} |\{c(u): u \in N(v)\}|+1,$$
where the minimum is taken over all proper vertex-colorings $c$ of $G$.
\end{defi}

Thus $\psi(G)$ is the minimum number of colors that must appear in the most
colorful closed neighborhood of a vertex in any proper coloring that may
involve an arbitrary 
number of colors. It was shown in \cite{EFHKRS} that there exist graphs $G$
with $\psi(G)=3$ and $\chi(G)>k$ for any positive integer $k$, where $\chi(G)$
denotes the chromatic number of $G$. 

Changing ``neighborhood'' to ``outneighborhood''
in the previous definition we arrive at the directed local chromatic number
(of a digraph) introduced in \cite{KPS}. For a directed graph $F$ let the
set of outneighbors of a vertex $v$ be $N_+(v)=\{u\in V(F): (v,u)\in E(F)\}$. 
By a proper vertex-coloring of a directed
graph we mean a proper vertex-coloring of the underlying undirected graph. 

\begin{defi} {\rm (\cite{KPS})}
The directed local chromatic number of a directed graph $F$ is defined as 
$$\psi_d(F)=\min_c\max_{v\in V(F)}\{c(u): u\in N_+(v)\}+1,$$
where $c$ runs over all proper vertex-colorings of $F$. 
\end{defi}

The directed local chromatic number of a digraph is always less than or equal
to the local chromatic number of the underlying undirected graph and we
obviously have equality if our digraph is symmetrically directed, i.e., for
every ordered pair $(u,v)$ of the vertices $(u,v)$ is an edge if and only if
$(v,u)$ is an edge. A digraph $F=(V,E)$ is called {\em oriented} if the
contrary is true: $(u,v)\in E$ implies $(v,u)\notin E$. An {\em orientation}
of an undirected graph $G$ is an oriented graph $\hat G$ that has $G$ as its
underlying undirected graph.

It is a natural question whether every undirected graph $G$ has an orientation 
the directed local chromatic number of which achieves the local chromatic
number of $G$. Currently we know very little about this question. (See
\cite{SaSi} for a problem of similar flavor: the relation of 
Shannon capacity and the maximum possible Sperner capacity of its
orientations.) 

In this paper we explore the other extreme: what is the minimum possible
directed local chromatic number that an orientation of a graph can attain. 

In the following section we give some more definitions and summarize some 
facts about the investigated parameters. In Section~\ref{shift} we investigate
shift graphs. We observe that they have an orientation with directed local
chromatic number $2$ and prove that their local chromatic number can be
arbitrarily large, in particular, it differs at most $1$ from their chromatic
number. We also consider the behavior of a symmetrized variant of shift
graphs.  

In section~\ref{blike} we
concentrate on {\em Borsuk-like graphs}: these are
graphs the chromatic number of which can be determined by applying Lov\'asz's
topological method (cf.\ \cite{LLKn}), while, at the same time they admit
optimal colorings where no short odd length walks exist that start and end in
the same color class. Several graphs have this property. In \cite{ST} we have
shown that the local chromatic number of these graphs is around one half of
their chromatic number. Here we show that the minimum directed
local chromatic number of a Borsuk-like graph of appropriate parameters is
about one quarter of its chromatic number.

\section{Minimum and maximum directed local chromatic number}

It is natural to define the following extreme values of $\psi_d(G)$.

\begin{defi}
For an undirected graph $G$ we define the {\em minimum directed local
  chromatic number} as 
$$\psi_{\rm {d,min}}(G):=\min_{\hat G}\psi_d(\hat G)$$
and the {\em maximum directed local chromatic number} as 
$$\psi_{\rm {d,max}}(G):=\max_{\hat G}\psi_d(\hat G),$$
where $\hat G$, in both cases, runs over all orientations of $G$. 
\end{defi}

It is obvious that $\psi_{\rm {d,max}}(G)\leq\psi(G)$. Equality holds for
complete 
graphs (by the transitive orientation), and more generally, for all graphs with
equal chromatic and clique number, thus for all perfect graphs, in
particular. A less obvious example for equality is given by Mycielski
graphs, see Proposition~\ref{prop:Myc} in Section~\ref{blike}.
We do not know whether equality holds for all graphs. 

Our main concern here, however, is the behavior of
$\psi_{\rm {d,min}}(G)$. Clearly, if the graph has any edge, then $\psi_{\rm
  {d,min}}(G)$ is already at least $2$. We will see in the next section that
there are graphs with $\psi_{\rm {d,min}}(G)=2$ and $\psi(G)$ arbitrarily
large.  

\medskip
To conclude this section we give an easy estimate on $\psi_{\rm d,min}(G)$ in
terms of $\chi(G)$. 
Recall that a homomorphism from graph $G$ to another graph $H$ is a mapping
$f:V(G)\to V(H)$ 
such that whenever $\{a,b\}$ is an edge of $G$, then $\{f(a),f(b)\}$ is an
edge of $H$. For a detailed introduction to graph homomorphisms, see
\cite{HN}.

The following relation of $\psi_{\rm {d,min}}$ to the chromatic number is
immediate.  

\begin{prop}\label{compl}
For every graph $G$
$$\psi_{\rm {d,min}}(G)\leq \left\lfloor\chi(G)\over2\right\rfloor+1.$$
If $G$ has equal clique number and chromatic number, then equality holds. 
\end{prop}

\proof
Let $G$ be a graph with chromatic number $r$, which means that there is a
homomorphism from $G$ to $K_r$. Orient the edges of $K_r$ so that the maximum
outdegree become as small as possible. Clearly, this minimal maximum outdegree
is $\lfloor r/2\rfloor$.

Let $c:V(G)\to V(K_r)$ be an optimal coloring of $G$. For each edge $\{u,v\}$
of $G$ orient it from $u$ to $v$ if and only if the edge $\{c(u),c(v)\}$ of
$K_r$ is oriented from $c(u)$ to $c(v)$ above. The set of colors in the
outneighborhood of a vertex $v$ of $G$ will be the set of outneighbors of
$c(v)$ in $K_r$. This proves $\psi_{\rm {d,min}}(G)\leq \lfloor {r\over
2}\rfloor+1.$ 

If the clique number of $G$ is also $r$ then some vertex of an $r$-clique of
$G$ must have at least $\lfloor{r\over 2}\rfloor$ other vertices of this
clique in its outneighborhood. Since all these must have different colors,
$\psi_{\rm {d,min}}(G)\ge \lfloor {r\over 2}\rfloor+1$ in this case.   
\qed

\section{Shift graphs}\label{shift}

Shift graphs were introduced by Erd\H{o}s and Hajnal \cite{EHsh}. 

\begin{defi} \label{defi:sh1} {\rm (\cite{EHsh})}
The shift graph $H_m$ is defined on the ordered pairs $(i,j)$ satisfying
$1\leq i<j\leq m$ as vertices and two pairs $(i,j)$ and $(k,\ell)$ form an
edge if and only if $j=k$ or $\ell=i$.
\end{defi}

Note that $H_m$ is isomorphic to the line graph of the transitive tournament
on $m$ vertices. It is well-known (see,
e.g., \cite{LLpr}, Problem 9.26) that $\chi(H_m)=\lceil\log_2 m\rceil$.  

Shift graphs are relevant for us for two different reasons. One is what we 
already mentioned in the Introduction that their minimum directed local
chromatic number is much below their local chromatic number. The other reason
is explained below. 

While the local chromatic number is obviously bounded from above by the
chromatic number, in 
\cite{KPS} it was shown to be bounded from below by the fractional chromatic
number. This motivated the study of the local chromatic number for graphs
with a large difference between the latter two bounds (see
\cite{ST}). 
Determining the chromatic number of such graphs often requires special tricks
as one needs some lower bound that is not a lower bound for the fractional
chromatic number. In case of Kneser graphs this difficulty was overcome by
Lov\'asz \cite{LLKn} thereby introducing his topological method that was later
successfully applied also for other graph families with the above
property. Examples 
include Schrijver graphs (\cite{Schr}) and generalized Mycielski graphs
(\cite{Stieb, GyJS}). See also \cite{Mat} for an excellent introduction to
this method. 

In \cite{ST} (see also \cite{STV}) we investigated the local chromatic number
of graphs for which the chromatic number is far from the
fractional 
chromatic number and can be determined by a particular implementation of the
topological method. If this implementation gave $t$ as a lower bound of the
chromatic number, we called a graph {\em topologically $t$-chromatic}, and
showed 
that if a graph is topologically $t$-chromatic, then $\lceil t/2\rceil+1$ is
an often tight lower bound for its local chromatic number. 

For shift graphs this topological lower bound for the chromatic number is not
tight 
(except for some very small initial cases), in other words they are not
topologically $t$-chromatic for $t$ being the actual chromatic number, see
Proposition~\ref{notop} below. On the other hand, shift graphs do have the
property that there is a large gap between their fractional and ordinary
chromatic numbers. Thus the above mentioned result of \cite{KPS} equally
motivates the 
investigation of their local chromatic number while the methods of \cite{ST,
  STV} cannot give good bounds for it. 

To see that the fractional chromatic number $\chi_f(H_m)$ is small it is
worth defining the {\it symmetric shift graph} $S_m$ that contains all ordered
pairs $(i,j)$ where $1\leq i, j\leq m$, $i\ne j$, as vertices (i.e., $(i,j)$ is
a vertex 
even if $i>j$) and $(i,j)$ and $(k,\ell)$ are adjacent again if $j=k$ or
$\ell=i$. (Note that $S_m$ is the line graph of the complete directed graph on
$m$ vertices.) It is obvious that $S_m$ is vertex-transitive, thus
$\chi_f(S_m)={|V(S_m)|\over \alpha(S_m)}$ (cf., e.g. \cite{SchU}), where
$\alpha(G)$ stands for the independence number of graph $G$. Since
  $\alpha(S_m)=\lceil{m\over 2}\rceil\lfloor{m\over 2}\rfloor$ (vertices
  $(i,j)$ with $i\leq\lceil{m\over 2}\rceil<j$ form an independent set of this
  size and one easily sees that no larger one can be formed), we get
  $\chi_f(H_m)\leq \chi_f(S_m)={{m(m-1)}\over {\lceil{m\over
        2}\rceil\lfloor{m\over 2}\rfloor}}<4$, where the
  first inequality follows from $H_m$ being a subgraph of $S_m$. 

Thus by the inequalities $\chi_f(H_m)\leq\psi(H_m)\leq\chi(H_m)$ the value of
$\psi(H_m)$ could be anywhere between $4$ and $\lceil\log_2 m\rceil$. Now we
show that the lower bound cannot be improved by the methods used in \cite{ST}. 

The lower bound on $\psi(G)$ in \cite{ST} mentioned above is
proven by showing (cf. also \cite{kyfan2} for a special case), 
that if $G$ is a topologically $t$-chromatic graph, then whatever way we color
its vertices properly (with any number of colors, thus the coloring need not
be optimal) there always appears a complete bipartite subgraph $K_{\lceil
  t/2\rceil, \lfloor t/2\rfloor}$, all $t$ vertices of which get a different
color. Though we do not give here the exact definition of topological
$t$-chromaticity, it makes sense to state the following proposition that can
be proven using the result just described. We remark that topological
$t$-chromaticity is a monotone property, that is, it implies topological
$(t-1)$-chromaticity. 

\begin{prop} \label{notop}
The graph $H_m$ is not topologically 4-chromatic and $S_m$ is not
topologically $5$-chromatic.
\end{prop}

\proof
Let us color the vertex $(i,j)$ with color $i$. This gives a proper
coloring of $H_m$.

One can easily check that if two vertices $(i_1,j_1)$ and $(i_2,j_2)$ of $H_m$
have two common neighbors $(k_1,\ell_1)$, $(k_2,\ell_2)$, then either
$j_1=j_2=k_1=k_2$ or $i_1=i_2=\ell_1=\ell_2$. Thus $H_m$ can be properly
colored in such a way it has no $K_{2,2}$ subgraph with all four vertices
receiving a different color. By the above described result in \cite{ST}, this
implies that $H_m$ is not topologically $4$-chromatic.

The same coloring (assigning color $i$ to the vertex $(i,j)$) is also a proper
coloring of $S_m$ but here for $m\ge4$ some $K_{2,2}$ subgraphs (like the one
consisting of the vertices $(1,2)$, $(2,3)$, $(3,4)$, $(4,1)$) receive four
distinct colors. However no $K_{2,3}$ subgraphs receive five distinct colors,
so by the same quoted result $S_m$ is not topologically $5$-chromatic.
\qed

We remark that $S_m$ is not even
topologically $4$-chromatic, but to see this is beyond the scope of the
present paper because every proper coloring of $S_4$ makes a $K_{2,2}$
subgraph (a 4-cycle) receive four distinct colors. Every non-bipartite graph
is topologically $3$-chromatic, so the graphs $H_m$ for $m\ge5$ and $S_m$ for
$m\ge3$ are topologically $3$-chromatic.

Although the local chromatic number of shift graphs could be as low as $3$ if 
considering only the topological lower bound of the local chromatic number
given in \cite{ST}, the main result of this section below states
that it is much higher. 

\begin{thm} \label{log}
We have 
$$\psi(H_m)
=\chi(H_m)$$ whenever $2^k+2^{k-1}<m\leq
2^{k+1}$ for some positive integer $k$. If $2^k<m\leq 2^k+2^{k-1}$ holds for
some $k$ instead, then we have 
$$
\chi(H_m)-1\leq \psi(H_m)\leq \chi(H_m)
.$$ 
\end{thm}

We prove this theorem in Subsection~\ref{proof of Th3}. 
It shows not only that the local chromatic number of shift graphs is close to
their chromatic number but also that the gap between the directed local
chromatic number and the local chromatic number of the underlying undirected
graph can be arbitrarily large. This statement follows when comparing
Theorem~\ref{log} to the following simple observation. (For the appearance of
more general shift graphs in a similar context, see the starting example in
\cite{EFHKRS}.) 

\begin{prop}\label{shmind}
$$\psi_{\rm d,min}(S_m)=\psi_{\rm d,min}(H_m)=2.$$
\end{prop}

\proof
As $H_m$ is a subgraph of $S_m$ and $\psi_{\rm d,min}(H_m)\ge 2$ is obvious,
it is enough to prove $\psi_{\rm d,min}(S_m)\leq 2$. Let $\tilde{S}_m$ be the
oriented version of $S_m$ in which edge $\{(a,b),(b,c)\}$ is oriented from
vertex $(a,b)$ to vertex $(b,c)$ whenever $a$, $b$ and $c$ are distinct while
we choose arbitrarily when orienting the edge between the vertices $(a,b)$ and
$(b,a)$ for $a\ne 
b$. Color each vertex $(x,y)$ by its first element $x$. Let $(a,b)$ be an
arbitrary vertex and observe that every element of its outneighborhood is
given color $b$. This shows $\psi_d(\tilde{S}_m)\leq 2$ thereby proving the
statement. 
\qed

Note the easy fact, that if we modify the directed graph $\tilde{S}_m$ in the
above proof so that for edges $\{(a,b),(b,a)\}$ we include both orientations
then the so obtained graph $\hat{S}_m$ is a homomorphism universal graph: it
has the property that a digraph $F$ admits a coloring with $m$ colors attaining
$\psi_d(F)\leq 2$ if and only if there exists a homomorphism from $F$ to
$\hat{S}_m$. (With the notation of \cite{KPS} $\hat{S}_m$ is just the graph
$U_d(m,2)$.) We will refer to the graphs $\hat{S}_m$ as the {\it symmetric
  directed shift graph}s.

\subsection{Bollob\'as-type inequalities}

A key observation in proving Theorem~\ref{log} will be the close connection
between local colorings of shift graphs and cross-intersecting set
systems. Here we state two classical results about the latter that will be
relevant for us. The first of these is due to Bollob\'as.

\begin{thm} \label{Bollin} {\rm (\cite{Boll})}
Let $A_1,\dots, A_m$ and $B_1,\dots, B_m$ be finite sets 
satisfying the property that $A_i\cap B_j\neq\emptyset$ for all $1\leq i,
j\leq m$ with $i\ne j$, while $A_i\cap B_i=\emptyset$ for all $1\leq i\leq
m$. Then
$$\sum_{i=1}^m {{|A_i|+|B_i|}\choose |A_i|}^{-1}\leq 1.$$
\end{thm}

Note that if $|A_i|=r$ and $|B_i|=s$ holds for all $i$ then the above
statement implies $m\leq {{r+s}\choose r}$. This consequence is generalized by
Frankl as follows. 

\begin{thm} \label{Frin} {\rm (\cite{PF})}
Let $A_1,\dots, A_m$ and $B_1,\dots, B_m$ be sets satisfying
$|A_i|=r, |B_i|=s, A_i\cap B_i=\emptyset$ for all $1\leq i\leq m$, and the
additional property that $A_i\cap B_j\neq\emptyset$ whenever $1\leq i<j\leq
m$. 
Then $$m\leq {{r+s}\choose r}.$$ 
\end{thm}

We remark that further relaxing the condition 
$A_i\cap B_j\neq\emptyset$ whenever $1\leq i<j\leq m$ to 
$1\leq i<j\leq m \Rightarrow (A_i\cap B_j\neq\emptyset\ {\rm or}\ A_j\cap
B_i\neq\emptyset),$ 
we arrive to a problem that, by our current knowledge, is not completely
solved for $r,s\ge 2$, cf. \cite{Tuza}. 

\medskip

The following lemma shows the connection between our problem and 
cross-intersecting set systems. 

\begin{lem} \label{Beq}
The inequality $\psi(H_m)\leq k$ is equivalent to the following statement. 
There exist finite sets, $A_1,\dots,A_m$ and $B_1,\dots, B_m$
such that $A_i\cap B_i=\emptyset$ for all $1\leq i\leq m$ and
for all $1\leq i<j\leq m$ we have $A_i\cap B_j\neq\emptyset$ and
$|A_j\cup B_i|\leq k-1$. 
\end{lem}

\proof
Assume first that $\psi(H_m)\leq k$ and let $c: V(H_m)\to \mathbb N$ be a
proper coloring that attains the local chromatic number. 
For each $1\leq i\leq m$ form the sets $A_i, B_i$
by $A_i:=\{c(i,\ell):i<\ell\le m\}, B_i:=\{c(\ell,i): 1\le\ell<i\}.$ 
Since the coloring is proper we must have
$A_i\cap B_i=\emptyset$ for all $i$. For $1\leq i<j\leq m$ we have $c(i,j)\in
A_i\cap B_j$, thus we have $A_i\cap B_j\neq\emptyset$ for all $i<j$.
A given vertex $(i,j)$ of $H_m$ is adjacent to the vertices $(\ell,i)$ and
$(j,q)$ where $\ell<i<j<q$. By our condition on the local chromatic number
this implies $|B_i\cup A_j|\leq k-1$ for all $i<j$. 

On the other hand, if $A_1,\dots, A_m, B_1,\dots, B_m$ with the above
properties exist, then we can define the coloring $c$ of the vertices of $H_m$
as follows.  
For each vertex $(i,j)\in V(H_m)$ let $c(i,j)$ be an arbitrary element of the
nonempty set $A_i\cap B_j$. As $A_i\cap B_i=\emptyset$ for all $i$ this
coloring is proper. By $|A_j\cup B_i|\leq k-1$ the local chromatic number
attained by this coloring is at most $k$. 
\qed

\subsection{Proof of Theorem~\ref{log}}\label{proof of Th3}

We will show that if the sets $A_1,\dots,A_m$ and
$B_1,\dots, B_m$ satisfy the conditions in 
Lemma~\ref{Beq}, then $m\leq 2^k+2^{k-1}$. By Lemma~\ref{Beq} and
$\chi(H_m)=\lceil\log_2 m\rceil$, this implies the
statement of Theorem~\ref{log}.

For obtaining the above upper bound on $m$ we partition the pairs 
$(A_i,B_i)$ according to the sizes of the sets $A_i,
B_i$. For every $0\leq r$ set  
$${\cal D}_1^{(r)}=\{i: 1\leq i\leq m, |A_i|=r, |A_i|+|B_i|<k\}$$
and 
$${\cal D}_2^{(r)}=\{i: 1\leq i\leq m, |A_i|=r, |A_i|+|B_i|\ge k\}.$$
Note that by its definition ${\cal D}_1^{(r)}=\emptyset$ for $r\ge k$ and
$|A_j\cup B_1|\le k-1$ for $1<j\le m$ implies $\cup_{r\ge k}{\cal
  D}_2^{(r)}\subseteq\{1\}$. 

Fix some $r\ge 0$. Notice that for each $i\in {\cal D}_1^{(r)}$ we have
$|B_i|\le k-1-r$ and add $k-1-r-|B_i|$ new elements to the set $B_i$ that do
not appear elsewhere. Denote the resulting set by $B_i'$. Note that the pairs
$(A_i,B_i')$ for $i\in{\cal D}_1^{(r)}$ satisfy the conditions in Frankl's
Theorem~\ref{Frin} (with $s=k-1-r$), implying $|{\cal D}_1^{(r)}|\leq 
{k-1\choose r}$. This further implies
$$\sum_{r\ge0}|{\cal D}_1^{(r)}|\leq 2^{k-1}.$$

For bounding the size of sets ${\cal D}_2^{(r)}$ observe that the condition
$|A_j\cup B_i|\leq k-1$ satisfied for all $i<j$ is equivalent to $|A_j\cap
B_i|\ge |A_j|+|B_i|-k+1$. Fix some $0\le r<k$ and notice that for $i\in
{\cal D}_2^{(r)}$ we have $|B_i|\ge k-r$. Let $B_i'$ be an arbitrary subset of
$B_i$ of size $k-r$. The pairs $(A_i,B_i')$ for $i\in{\cal D}_2^{(r)}$ still
satisfy that $A_j\cap B_i'\neq\emptyset$ whenever $j>i$, while $A_i\cap
B_i'=\emptyset$ is also true. Thus the conditions of Theorem~\ref{Frin} hold
again (now with $s=k-r$ and by reversing the order of indices) implying
$|{\cal D}_2^{(r)}|\leq {k\choose r}.$
This further implies
$$\sum_{r\ge 0}|{\cal D}_2^{(r)}|
\leq \sum_{r=0}^{k-1}|{\cal D}_2^{(r)}|+1\leq 2^k.$$ 

\medskip
Thus we obtained $m=\sum_{r\ge0}|{\cal D}_1^{(r)}|+\sum_{r\ge0}|{\cal
D}_2^{(r)}|\leq2^k+2^{k-1}$ completing the proof. 
\qed

\subsection{Symmetric shift graphs}

In view of the above it is natural to ask what is the local chromatic number
of the symmetric shift graph $S_m$. We trivially have
$\psi(S_m)\ge\psi(H_m)$. In view of Theorem~\ref{log} this shows that
$\psi(S_m)$ is close to $\chi(S_m)=\min\left\{k:
{k\choose {\lceil k/2\rceil}}\ge m\right\}$ (see, e.g. \cite{LLpr}, Problem
9.26.), but this trivial observation allows for
an unbounded difference of the order $\log(\chi(S_m))$ or $\log\log m$. In view
of Theorem~\ref{log} it seems very unlikely that there could be
such a large gap between $\psi(S_m)$ and $\chi(S_m)$. In fact, we are inclined
to believe that both $\psi(S_m)$ and $\psi(H_m)$ coincides with the
corresponding chromatic numbers, $\chi(S_m)$ and $\chi(H_m)$, respectively.

In this subsection we apply the method of the preceding section to
improve the above trivial lower bound on $\psi(S_m)$. The improvement we
obtain is rather modest: we increase the lower bound by $1$ for some $m$.

The analogue of Lemma~\ref{Beq} is the following.

\begin{lem}\label{sBeq}
The inequality $\psi(S_m)\leq k$ is equivalent to the following statement. 
There exist finite sets $A_1,\dots,A_m$ and $B_1,\dots, B_m$
such that $A_i\cap B_i=\emptyset$ for all $1\leq i\leq m$ and
for all $1\leq i,j\leq m$ with $i\ne j$ we have $A_i\cap B_j\neq\emptyset$ and
$|A_i\cup B_j|\leq k-1$.
\end{lem} 

The proof is essentially identical to that of Lemma~\ref{Beq}, therefore we
omit it.

\begin{thm}\label{slog}
The local chromatic number of the symmetric shift graph $S_m$ satisfies
$$\psi(S_m)\ge \lceil \log_2 (m+2)\rceil.$$
\end{thm}

\medskip

\proof
We do the same as in the proof of Theorem~\ref{log}. By Lemma~\ref{sBeq} it is
enough to show that if $A_1,\dots,A_m$ and
$B_1,\dots, B_m$ are two families of finite sets satisfying the conditions
there, then $m\leq 2^k-2$. 

To this end we define ${\cal D}^{(r)}=\{i: 1\leq i\leq m, |A_i|=r\}$.

Note that for $r\ge k$ \ ${\cal D}^{(r)}=\emptyset$ follows from the 
condition $|A_i\cup B_j|\leq k-1$ for $i\neq j$. Similarly, ${\cal
D}^{(0)}=\emptyset$ follows from $A_i\cap B_j\neq\emptyset$ for $i\neq j$. 

Fix some $0<r<k$ and consider $i\in {\cal D}^{(r)}$. If $|B_i|>k-r$ let
$B_i'$ be an arbitrary subset of $B_i$ of size $k-r$, otherwise let
$B_i'=B_i$. The conditions imply that the pairs $(A_i,B_i')$ for $i\in{\cal
D}^{(r)}$ satisfy the conditions of Theorem~\ref{Bollin}. Since we have
$|A_i|=r$, $|B_i'|\leq k-r$ for all $i\in{\cal D}^{(r)}$, this further
implies $|{\cal D}^{(r)}|\le {k\choose r}$. Summing for
all $r$ we obtain 
$$m=\sum_{r=1}^{k-1}|{\cal D}^{(r)}|\leq 2^k-2$$
completing the proof. 
\qed

\subsection{A homomorphism duality result}

In this subsection we prove that the following homomorphism duality statement
(see 
\cite{HN} for more on this term) holds for symmetric directed shift graphs
$\hat S_m$ (see their definition after Proposition~\ref{shmind}). 
We need the notion of an alternating odd cycle, which is an oriented odd cycle
with exactly one vertex of outdegree one.
It was observed in \cite{KPS} that a directed odd cycle has directed local
chromatic number $3$ if and only if it contains an alternating odd cycle as a
subgraph. The following is a straightforward extension of this observation.

\begin{prop}
A directed graph $\hat G$ admits a homomorphism into $\hat S_m$ for some $m$
if and only if no alternating odd cycle admits a homomorphism to $\hat G$. 
\end{prop}

\proof
It is clear (and also contained in \cite{KPS}) that alternating odd cycles
have directed local chromatic number $3$. By the remark following the proof of 
Proposition~\ref{shmind} this implies that
there is no homomorphism from any alternating odd cycle to $\hat S_m$ for any
$m$, or to any graph that admits a homomorphism to a symmetric directed shift
graph $\hat S_m$ for some $m$. 

On the other hand, we claim that if $\psi_d(\hat G)>2$ (which is equivalent to 
$\hat G$ not having a homomorphism to any $\hat S_m$), then an alternating odd
cycle has a homomorphism to $\hat G$. (We remark that this also implies that
$\hat G$ contains an alternating odd cycle as a subgraph.) 
Indeed, call two vertices $u$ and $v$ {\em related}
if they both belong to the outneighborhood of the same vertex $w$. 
The transitive closure of this relation defines equivalence classes of the
vertices. Let us color the vertices according to the equivalence class they
belong to. Clearly, the outneighborhood of any vertex is monochromatic, so
$\psi_d(\hat G)>2$ implies that this is not a proper coloring of $\hat G$.
Let $a$ and $b$ be adjacent vertices in an equivalence class. 
There must be a sequence $a=u_0, u_1,\dots, u_h=b$ of vertices such
that $u_i$ is related to $u_{i+1}$ for $0\le i<h$. Let $w_i$ be the vertex
having both $u_i$ and $u_{i+1}$ in its outneighborhood. The vertices of an
alternating odd cycle of length $2h+1$ can be homomorphically mapped to $u_0,
w_0, u_1,w_1,\dots, u_h$ in this order. \qed

\section{Borsuk-like graphs}\label{blike}

Borsuk-graphs were also introduced by Erd\H{o}s and Hajnal \cite{EH}. 

\begin{defi}{\rm (\cite{EH})}
The Borsuk graph $B(n,\alpha)$ is defined for every positive integer $n$ and
$0<\alpha<2$ on the unit sphere $\mathbb S^{n-1}$ of the $n$-dimensional
Euclidean space as 
vertex set. Two vertices form an edge if their Euclidean distance is larger
than $\alpha$.
\end{defi}

It is easy to see that the statement $\chi(B(n,\alpha))\ge n+1$ is
equivalent with the celebrated Borsuk-Ulam theorem, see \cite{EH, LLBor}. It
is also well-known and easy to see, that if $\alpha$ is larger than a certain
threshold, than $n+1$ colors suffice: inscribe a regular simplex into $\mathbb
S^{n-1}$ and color each point of the sphere with the side of the simplex
intersected by the line segment joining this point to the center of the
sphere. Note that besides being proper this coloring has a further remarkable
property: for every $s\in\mathbb N$ there exists $\alpha_{n,s}<2$ such that if
$\alpha>\alpha_{n,s}$
then there is no walk of length $2s-1$ in $B(n,\alpha)$ between any pair of
vertices that have the same color. Several other interesting graphs also have
optimal colorings with this property, see \cite{ST}. 

\begin{defi}{\rm (\cite{ST}, cf. also \cite{BaumStieb})}
Let $s$ be a positive integer. 
A coloring $c$ of a graph $G$ is called {\em $s$-wide} if there is no walk of
length $2s-1$ in $G$ between any two vertices $u$ and $v$ with $c(u)=c(v)$.   
\end{defi}

Observe that $1$-wide colorings are exactly the proper colorings, while being
$2$-wide means that the neighborhood of each color class is
independent. Graphs with colorings of the latter property were investigated in
\cite{GyJS}. $3$-wide colorings were simply called {\em wide} in \cite{ST} as
they had a key role there in bounding the local chromatic number from
above. Namely, 
we proved in \cite{ST} that if a graph $G$ has a $3$-wide coloring with $t$
colors then $\psi(G)\leq \lfloor t/2\rfloor +2$. (To see that this bound is
sharp for several graphs, cf. \cite{ST, STV}.)  

\medskip

Recall that the Kneser graph ${\rm KG}(n,k)$ is defined for $n\ge2k$ on all
$k$-element subsets of the $n$ element set $[n]=\{1,\dots,n\}$ as vertex set
and two such subsets form an edge if they are disjoint. Their chromatic number
is $n-2k+2$ as conjectured by Kneser \cite{Kne} and proved by Lov\'asz
\cite{LLKn}. Schrijver found a very nice family of induced subgraphs of Kneser
graphs. They have the same chromatic number as the corresponding Kneser graphs
but at the same time they are also vertex color-critical. 

\begin{defi}\label{SGnk} {\rm (\cite{Schr})}
The {\em Schrijver graph} $\SG(n,k)$ is defined for $n\ge2k$ as
follows.
\begin{eqnarray*}
V(\SG(n,k))&=&\{A\subseteq [n]: |A|=k,\forall i:\ \{i,i+1\}\nsubseteq
A\ \ \hbox{\rm and}\ \ \{1,n\}\nsubseteq A\}\\
E(\SG(n,k))&=&\{\{A,B\}: A\cap B=\emptyset\}
\end{eqnarray*}
\end{defi}

The following generalization of Mycielski's construction \cite{Myc} appears
in several papers, see, e.g., \cite{GyJS, Stieb, Tar} for their chromatic
properties.

\begin{defi}\label{genMyc}
For a graph $G$ and integer $r\ge 1$ the generalized Mycielskian $M_r(G)$ of
$G$ is the graph on vertex set 
$$V(M_r(G))=\{(i,v): v\in V(G), 0\leq i\leq r-1\}\cup\{z\}$$
with edge set 
$$E(M_r(G))=\{\{(i,u),(j,v)\}: \{u,v\}\in E(G)\ {\rm and}\ i=j=0\ {\rm or}\
0\leq i=j-1\leq r-2\}\cup$$
$$\{\{(r-1,u),z\}: u\in V(G)\}.$$
\end{defi}

\smallskip

The Mycielskian $M(G)$ of a graph is identical to $M_2(G)$. The main property
of this construction is that while it does not change the clique number for
$r\ge 2$, the chromatic number of $M(G)$ is $1$ more than that of
$G$. We have $\chi(M_r(G))\le\chi(G)+1$ for an arbitrary $r$, but
$\chi(M_r(G))=\chi(G)$ can happen for $r\ge 3$ (an example is
$G=\bar C_7$, see
\cite{Tar}, or see \cite{Cs} for another example with fewer edges).
Stiebitz \cite{Stieb} proved, however, that Lov\'asz's topological lower
bound on the chromatic number is always $1$ more for $M_r(G)$ than for
$G$. Thus, if this bound is tight for $G$ then the chromatic number of $M_r(G)$
is $1$ larger than $\chi(G)$. Moreover, in this case this new bound is also
tight for $M_r(G)$, so this argument can be used recursively.

\medskip

The chromatic number of all the above graphs were determined by using the
topological method, in particular, the Borsuk-Ulam theorem, for getting the
appropriate lower bound, see \cite{LLKn, Schr, Stieb, GyJS} and also
\cite{Mat}.  
Another similarity between Schrijver graphs and generalized Mycielski 
graphs is that for any given chromatic number $\chi$ and parameter $s$ one can
find a member of either family with chromatic number $\chi$ having an $s$-wide
$\chi$-coloring. (We note that a topological
similarity of Schrijver graphs and their iterated generalized
Mycielskians that is not shared by Kneser graphs is that their so-called 
neighborhood complex, cf.\ \cite{LLKn, Mat}, is homotopy equivalent to a
sphere, see \cite{BjdeL, Stieb}.)

We conclude the introductory part of this section by stating a result about
the maximum directed local chromatic number of Mycielski graphs. It is a
rather straightforward generalization of Proposition 10 in \cite{ST}. Though
its proof is almost identical to that of this quoted result, we
include it for the sake of completeness. 

\begin{prop}\label{prop:Myc}
For any graph $G$ we have 
$$\psi_{\rm d,max}(M(G))\ge\psi_{\rm d,max}(G)+1.$$
In particular, if $\psi_{\rm d,max}(G)=\chi(G)$, then $\psi_{\rm
  d,max}(M(G))=\psi_{\rm d,max}(G)+1=\chi(M(G)).$ 
\end{prop}

\proof
First we give the orientation. Fix an orientation of $G$ that attains
$\psi_{\rm d,max}(G)$ and orient the subgraph of $M(G)$ induced by the
vertices $(0,v)$ accordingly. Orient each edge of the form $\{(1,u),(0,v)\}$
consistently with the corresponding edge $\{(0,u),(0,v)\}$, i.e., so that
either both have 
its head or both have its tail at the vertex $(0,v)$. Finally, orient all edges
$\{(1,u),z\}$ towards $z$. 

Now consider an arbitrary proper coloring $c:V(M(G))\to \mathbb N$. For a
subset $U\subseteq V(M(G))$ let $c(U):=\{c(u):u\in U\}$.
Consider also the modified coloring $c'$ of $G$ defined by 
$$c'(x)=\left\{\begin{array}{lll}c(0,x)&&\hbox{if }
c(0,x)\ne c(z)\\
c(1,x)&&\hbox{otherwise.}\end{array}\right.$$

It follows from the construction that $c'$ is a proper coloring of $G$, which
does not use the color $c(z)$. 

By our orientation of $G$ there is some vertex $v$ of $G$ for which
$|c'(N_+(v))|\ge \psi_{\rm d,max}(G)-1$. 
(Note that $N_+(.)$ and $N_+(.,.)$ here refer to outneighborhoods in the
considered orientations of $G$ and
$M(G)$, respectively.) 
If there is no vertex $u\in N_+(v)$ for
which $c(0,u)\ne c'(u)$, then the color $c(z)$ does not appear in the
outneighborhood of $(0,v)$ in $M(G)$. In this case the set $c(N_+(1,v))$
contains all 
the colors in $c'(N_+(v))$ plus the additional color $c(z)$. If, however,
there is 
some $u\in N_+(v)$ for which $c(0,u)\ne c'(u)$, then we have
$c(0,u)=c(z)$. In this case the set $N_+(0,v)$ contains all the colors
appearing in $c'(N_+(v))$ and also the additional color $c(z)$ as the color of
$(0,u)$. In either case, some vertex has at least 
$\psi_{\rm d,max}(G)$ colors in its
outneighborhood, proving $\psi_{\rm d,max}(M(G))\ge\psi_{\rm d,max}(G)+1.$    

The second statement trivially follows from the first using the well-known fact
$\chi(M(G))=\chi(G)+1$ and the obvious inequalities $\psi_{\rm
  d,max}(G)\leq\psi(G)\leq\chi(G)$. \qed

\medskip

Note that Proposition~\ref{prop:Myc} implies that $\psi_{\rm
  d,max}(G)=\psi(G)$ holds whenever $G$ is a Mycielski graph, that is a graph
obtained from a single edge by repeated use of the Mycielski construction. 
We also remark that unlike the analogous inequality for $\chi(G)$ or $\psi(G)$
the inequality $\psi_{\rm d,max}(M(G))\leq\psi_{\rm d,max}(G)+1$ does not seem
to be obvious. Though we do not have a counterexample we are not
completely convinced about its validity.

\subsection{Lower bound by topological t-chromaticity}\label{lowb}

As we have already mentioned in Section~\ref{shift} we called a graph 
topologically $t$-chromatic in \cite{ST} if a particular implementation of the
topological method gave $t$ as a lower bound for its chromatic number. We also
mentioned there that a result in \cite{ST} shows (cf. also \cite{kyfan2}) that
in every proper coloring of a topologically $t$-chromatic graph a complete
bipartite subgraph $K_{\lceil t/2\rceil, \lfloor t/2\rfloor}$ occurs, all $t$
vertices of which get a different color. This result was used in \cite{ST} to
bound $\psi$ from below. In a similar manner it also gives a lower bound for
$\psi_{\rm {d,min}}$.  

\begin{thm}\label{negyed}
If $G$ is a topologically $t$-chromatic graph with $t\ge2$, then
$$\psi_{\rm{d,min}}(G)\ge \lceil t/4\rceil+1.$$
\end{thm}

\proof
Let $G$ be a topologically $t$-chromatic graph, $c$ its proper coloring and
$D$ its multicolored complete bipartite subgraph whose existence is guaranteed
by the result mentioned above. The number of edges in $D$ is $\lceil
t/2\rceil\lfloor t/2 \rfloor$ implying that for any orientation of $D$ its
average outdegree is  
$(1/t)\lceil t/2\rceil\lfloor t/2 \rfloor$ the upper integer part of which is
$\lceil t/4\rceil$. Since all vertices of $D$ receive different colors, its
maximum outdegree vertex have at least $\lceil t/4\rceil$ different colors in
its outneighborhood in any orientation. This proves that $\psi_{\rm
  {d,min}}\ge \lceil t/4\rceil+1.$
\qed

\subsection{Upper bound by wide colorability}

\subsubsection{Graphs with chromatic number at most six}

If a graph $G$ is at most $3$-chromatic (but not edgeless), then
Proposition~\ref{compl} implies that its minimum directed local chromatic
number $\psi_{\rm {d,min}}(G)=2$. Below we will show that the same conclusion
holds for $4$-chromatic graphs with $2$-wide $4$-colorings. The same method
will be used to prove the sharpness of our topological lower bound for certain
graphs of chromatic number at most $6$.   

\medskip

The following notations and lemmas will be useful. 
Given a coloring $c$ of a graph $G$ for each vertex $v\in V(G)$ let
$S_c(v)=\{c(u): \{u,v\}\in E(G)\}$ and $s_c(v)=|S_c(v)|$. That is, $s_c(v)$ is
the number of colors given to the neighbors of $v$.

\begin{lem}\label{paros}
If $c$ is a $2$-wide coloring and $u,v$ are adjacent vertices
of a graph $G$ then $S_c(u)\cap S_c(v)=\emptyset$. In particular, if $c$ uses
$t$ colors, then $s_c(u)+s_c(v)\leq t$. 
\end{lem}

\proof
Assume indirectly that $S_c(u)\cap S_c(v)\neq\emptyset$, i.e., $u$ has a
neighbor $x$ and $v$ has a neighbor $y$ with $c(x)=c(y)$. But then the walk
$xuvy$ connects vertices of the same color and contradicts the
assumption that $c$ is $2$-wide. This proves the first statement of the
lemma, that obviously implies the second one completing the proof.
\qed

\begin{lem}\label{hh}
If a graph $G$ has a $2$-wide coloring using $2h$ colors with $h\ge2$, then
$\psi_{\rm{d,min}}(G)\leq h$.
\end{lem}

\proof
Consider $G$ as colored by a fixed $2$-wide $2h$-coloring
$c:V(G)\to H$ with $|H|=2h$.

Let us consider the subgraph $G'$ obtained from $G$ by removing all vertices
$u\in V(G)$ with $s_c(u)<h$. We claim that $G'$ has an orientation $\hat
G'$ such that the outneighborhood of any vertex receives at most $\lceil
h/2\rceil$ distinct colors by $c$.

Indeed, by Lemma~\ref{paros} if $\{u,v\}$ is an
edge of $G'$, then $S_c(u)$ and $S_c(v)$ are complementary sets of colors,
each of size $h$. So each nontrivial component of $G'$ is a bipartite graph
with one side containing vertices $u$ with $S_c(u)=H_1$ for some fixed set
$H_1$ of $h$ colors and with the other side containing vertices $v$ with
$S_c(v)=H_2=H\setminus H_1$. Clearly, the vertices in the former side receive
colors in $H_2$, while vertices on the latter side have colors in $H_1$. To
prove the claim it is enough to find a suitable orientation for each of the
components 
separately, so let us fix $H_1$ and $H_2$. Consider the complete bipartite
graph $K_{H_1,H_2}$ on the vertex set $H$ consisting of the edges connecting
elements of $H_1$ and $H_2$. Orient the edges of this graph, so that every
outdegree is at 
most $\lceil h/2\rceil$. Now orient the edge $\{u,v\}$ in this connected
component of $G'$ according to the orientation of $\{c(u),c(v)\}$ in
$K_{H_1,H_2}$. Clearly, this orientation satisfies the requirement of the
claim.

Having found the orientation $\hat G'$, extend it to an orientation $\hat G$
of $G$ by orienting each edge of $G$ not in $G'$ away from a vertex $u$ with
$s_c(u)<h$. The outneighborhood of a vertex in $G'$ is the same in $\hat G$ and
in $\hat G'$, so it receives at most $\lceil h/2\rceil\le h-1$ colors at
$c$. For the rest of the vertices of $G$ their entire neighborhood receives at
most $h-1$ colors, so we have $\psi_{\rm d}(\hat G)\le h$. This completes the
proof of the lemma.

Notice that the coloring establishing the bound on the directed local
chromatic number is the $2$-wide coloring itself. \qed

\begin{cor}\label{4chrom}
If a non-edgeless graph $G$ has a $2$-wide $4$-coloring, then $\psi_{\rm
  {d,min}}(G)=2$.  
\end{cor}

\proof
The statement immediately follows by applying Lemma~\ref{hh} with $h=2$. 
\qed
  
\medskip

\begin{cor}
If a topologically $5$-chromatic graph $G$ has a $2$-wide coloring using at
most $6$ colors, then $\psi_{\rm {d,min}}(G)=3$.  
\end{cor}
 
\proof
Theorem~\ref{negyed} implies $\psi_{\rm {d,min}}(G)\ge 3$.  
Lemma~\ref{hh} implies $\psi_{\rm {d,min}}(G)\le 3$.
\qed

\subsubsection{General upper bound}

In this section we improve Lemma~\ref{hh} so that the upper bound it gives will
match the lower bound of Theorem~\ref{negyed} for several graphs of higher
(local) chromatic number. For this we need to assume the existence of $s$-wide
colorings for larger values of $s$. In \cite{ST} the minimal universal graphs
for $s$-wide $t$-colorability were found. (Cf. \cite{GyJS} for some larger
universal graphs for this property.) We will use them here.

\begin{defi} Let $s\ge1$ and $t\ge2$ be integers.
The vertex set of the graph $W(s,t)$ consists of the functions
$f:\{1,\ldots,t\}\to\{0,1,\ldots,s\}$ satisfying that $f(i)=0$ holds for
exactly one index $i$ and $f(i)=1$ holds for at least one index $i$. Two
vertices $f$ and $g$ are connected in $W(s,t)$ if for every $i$ one has
$|f(i)-g(i)|=1$ or $f(i)=g(i)=s$.

The natural coloring of $W(s,t)$ assigns the color $i$ to the vertex $f$ if
$f(i)=0$.
\end{defi}

\begin{lem}\label{uni}{\rm (\cite{ST}, cf. also \cite{BaumStieb})} For $s\ge1$
    and $t\ge2$ the natural 
coloring of $W(s,t)$ is an $s$-wide $t$-coloring. A graph $G$ admits an
$s$-wide $t$-coloring if and only if there is a homomorphism from $G$ to
$W(s,t)$.
\end{lem}

\begin{thm}\label{swide}
For every $t\in \mathbb N$ there is an $s=s_t$ for which the following is true.
If a graph $G$ has an $s$-wide coloring with $t$ colors then $\psi_{\rm
  {d,min}}(G)\leq \lceil{t/4}\rceil+1.$ 
\end{thm}

\medskip

\proof We will find an orientation $\hat W$ of $W(s,t)$ with directed local
chromatic number bounded above by $\lceil t/4\rceil+1$. This is enough by
Lemma~\ref{uni} and the trivial observation that if there is a homomorphism
from a graph $G$ to another graph $W$, then we have $\psi_{\rm
d,min}(G)\le\psi_{\rm d,min}(W)$.

Let $\chi$ stand for the natural coloring of $W(s,t)$. This is the coloring
establishing our bound on $\psi_{\rm d}(\hat W)$. We
write $\tau$ for
$\lceil t/4\rceil$. We will define a set
$S(f)$ of colors for every vertex $f$ of $W(s,t)$. We make sure that
\begin{enumerate}
\item$|S(f)|\le \tau$ for every vertex $f$ and
\item either $\chi(f)\in S(g)$
or $\chi(g)\in S(f)$ holds for every edge $\{f,g\}$ of $W(s,t)$.
\end{enumerate}
We obtain the
orientation $\hat G$ by orienting an edge from $f$ to $g$ only if $\chi(g)\in
S(f)$. Property~2 ensures that all edges of $W(s,t)$ can be oriented
this way. Property~1 makes sure that the natural coloring $\chi$
establishes $\psi_{\rm d}(\hat G)\le \tau+1$. So finding the sets $S(f)$ with
these properties completes the proof of the theorem.

Let us fix a vertex $f$ of $W(s,t)$. Let $c=\chi(f)$, $E=\{1\le i\le
t:f(i)\hbox{ is even}\}$ and $O=\{1\le i\le t:f(i)\hbox{ is odd}\}$.
For $1\le i\le t$ let $p_i=\sum_{j\in E, j\le i}(s-f(j))$ and $q_i=\sum_{j\in
O,j\le i}(s-f(j))$. Note that $f(c)=0$, so $p_t\ge s$ and as there is an index
$i$ with $f(i)=1$ we have $q_t\ge s-1$.

The idea is to represent the colors in $E$ and in $O$ as points of
the real interval $[0,1]$ and orient the edges from $f$
towards those other vertices whose color in the natural coloring is
represented by a point which is circularly (that is, when identifying $1$ with
$0$) ``somewhat to the right'' from the point representing  
the color of $f$. To make this
orientation consistent for the different vertices of $W(s,t)$ we apply
appropriate weightings to determine the distances between the points
representing different colors. These weights will depend on the actual values
$f(i)$ for each color $i$ that measure the length of the shortest walk in
$W(s,t)$ from $f$ to a vertex of color $i$ in the natural coloring. 

If $f(1)$ is even, we set $P_i=(p_i-(s-f(1))/2)/p_t$ and
$Q_i=q_i/q_t$ for $1\le i\le t$. If $f(1)$ is odd we set $P_i=p_i/p_t$ and
$Q_i=(q_i-(s-f(1))/2)/q_t$. We have $0\le P_i,Q_i\le1$.

Note that $s-f(1)$ is a summand in one of $p_i$
and $q_i$ and the correction term of subtracting half of this summand is a
technicality that we will need to be able to prove the theorem also in the
case when $t$ is divisible by $4$.

Let $\varepsilon=t/(s-1)$. Note that $\varepsilon>0$ can be made arbitrarily
close to zero by choosing $s$ large enough for a fixed $t$. We express this
relationship simply by saying $\varepsilon$ is {\em close} to zero and will
use this term in similar meaning later in this proof.

In case there are at most $\tau$ indices $i$ with $f(i)=1$ we define $S(f)$ to
be the set of these indices. Otherwise we compute $D_i=Q_i-P_c+2\varepsilon$
for all indices $i$ with $f(i)=1$ and let $S(f)$ be formed by the $\tau$ indices
that have the smallest fractional parts $X_i=D_i-\lfloor D_i\rfloor$.

Property 1 is clear from the definition. In the rest of this proof we
establish property 2 if $s$ is large enough.

Assume for a contradiction that the vertices $f$ and $f'$ are connected in
$W(s,t)$ but property 2 fails for this edge. Let $c$, $p_i$, $q_i$, $P_i$,
$Q_i$, $D_i$ and $X_i$ be the above defined values for the vertex $f$ and let
$c'$, $p_i'$, $q_i'$, $P_i'$, $Q_i'$, $D_i'$ and $X_i'$ be the
corresponding values for $f'$.

First observe that as $f$ and $f'$ are connected $|f(i)-f'(i)|\le1$ for all
$i$ while $f(i)$ and $f'(i)$ are of different parity unless $f(i)=f'(i)=s$. 
This shows that $|p_i-q_i'|\le t$ and $|q_i-p_i'|\le t$ for all $i$. Easy
calculation shows that with our lower bound on $p_t$ and $q_t$ this implies 
$|P_i-Q_i'|\le 2\varepsilon$ and similarly 
$|Q_i-P_i'|\le 2\varepsilon$.

We have $f(c)=0$, $f'(c)=1$, $f'(c')=0$ and $f(c')=1$. By the formula defining
$D_i$ we have $0\le
D_{c'}+D'_c\le8\varepsilon$. For the fractional parts this means
$X_{c'}+X'_c\le 1+8\varepsilon$. We assumed that property 2 is violated, so
there are $\tau$ indices $i$ with $f(i)=1$ and $X_i<X_{c'}$ and similarly, for
$\tau$ indices $j$ we have $f'(j)=1$ and $X'_j<X'_c$.

It is easy to see that the values $X_i$ for indices satisfying $f(i)=1$ are
separated from each other 
by at least $(s-1)/q_t$, so we have $X_{c'}\ge \tau(s-1)/q_t$ and
therefore $q_t\ge \tau(s-1)/X_{c'}$. Similarly we have $q'_t\ge
\tau(s-1)/X'_c$. Using also the bound on $X_{c'}+X'_c$ we obtain
$q_t+q'_t\ge4\tau(s-1)/(1+8\varepsilon)$.

Notice that no index $i$ can contribute to both $q_t$ and $q'_t$. This is
because either one of $f(i)$ or $f'(i)$ is even and thus does not contribute
or if $f(i)=f'(i)=s$ is odd, then both contributions are zero. Those indices
that do contribute to either $q_t$ or $q'_t$ contribute at most $s-1$, so we
have $q_t+q'_t\le t(s-1)$. If $t<4\tau$ and $\varepsilon$ is small enough this
contradicts our lower bound on $q_t+q'_t$ and thus completes the proof of
property 2 in the $t<4\tau$ case.

In the tight $t=4\tau$ case we have to work more for the contradiction. We
still 
have $t(s-1)\ge q_t+q'_t\ge4\tau(s-1)/(1+8\varepsilon)$, but this inequality
does not lead directly to a contradiction. Let $\alpha>0$. If $\varepsilon$ is
small enough (the threshold depends on $t$ and $\alpha$), then
it yields that $q_t+q'_t\ge(t-\alpha)(s-1)$ and therefore, since any index can
contribute at most $(s-1)$ to one of $q_t$ and $q'_t$, each index $i$ must
contribute at least $(1-\alpha)(s-1)$ to $q_t$ or $q'_t$ (in other words
$f(i)$ must be small relative to $s$). Also, from $t(s-1)\ge q_t+q'_t\ge
\tau(s-1)/X_{c'}+\tau(s-1)/X'_c$ one obtains $1/X_{c'}+1/X'_c\le 4$, thus
$X_{c'}$ must be close to $1/2$. (Recall that this means that fixing $t$
and choosing $s$ large enough $|X_{c'}-1/2|$ can be made arbitrarily
small.) Now from $q_t\ge \tau(s-1)/X_{c'}$ (and $s$ large enough) it follows
that at least $2\tau$ indices contribute
to $q_t$ and similarly, at least $2\tau$ indices contribute to $q'_t$, so by
$4\tau=t$, exactly $2\tau$ indices contribute to each. Thus exactly
$2\tau$ indices contribute to $p_t$, as well.

We can assume by symmetry that $f(1)$ is odd: otherwise switch the roles of
$f$ and $f'$. Now we can estimate $P_c$ and
$Q_{c'}$. We have $P_c=p_c/p_t$ and, by the above, this is close to $2k/t$,
where $k=|\{1\le i\le c:f(i)\hbox{ is even}\}|$. We have
$Q_{c'}=(q_{c'}-(s-f(1))/2)/q_t$, and, similarly, this is close to
$(2\ell-1)/t$, where 
$\ell=|\{1\le i\le c':f(i)\hbox{ is odd}\}|$. This makes
$D_{c'}=Q_{c'}-P_c+2\varepsilon$ close to $(2\ell-2k-1)/t$. Here the
numerator is odd, the denominator is the fixed value $t$ divisible by $4$, so
the fractional part $X_{c'}$ of this number cannot be close to $1/2$. This
provides the contradiction proving property 2 and completing the proof of the
theorem.
\qed

\medskip

In the following corollaries $s=s_t$ always refers to the $s_t$ of
Theorem~\ref{swide}. 

\begin{cor}\label{topswi}
If $G$ is a topologically $t$-chromatic graph that has an $s$-wide
$t$-coloring for the value $s=s_t$, then $\psi_{\rm
d,min}(G)=\lceil{t/4}\rceil+1.$
\end{cor}

\proof
Follows from Theorems~\ref{negyed} and \ref{swide}.
\qed

\medskip

Finally, we specify two interesting special cases of
Corollary~\ref{topswi}. They rely on the topological and wide colorability
properties of the relevant graphs established in \cite{ST}. 

\begin{cor}
If $t=n-2k+2$ is fixed and $n\ge (2s-2)t^2-(4s-5)t$ for $s=s_t$, then
$$\psi_{\rm {d, min}}(\SG(n,k))=\left\lceil {t\over 4}\right\rceil+1.$$  
\end{cor}

\proof
It is shown in Lemma 5.1 of \cite{ST} that if the conditions in the statement
are satisfied, then $\SG(n,k)$ admits an $s$-wide $t$-coloring. Thus the
statement is implied by Theorem~\ref{swide} and the fact that $\SG(n,k)$ is
topologically $t$-chromatic (cf. \cite{Mat, Schr} or Proposition 8 in
\cite{ST}). \qed

\begin{cor}
If $G$ is a topologically $t$-chromatic graph admitting an $s$-wide
$t$-coloring for $s=s_t$ 
and $r\ge3s-2$, then
$$\psi_{\rm {d, min}}(M_r(G))=\left\lceil{t+1\over 4}\right\rceil+1.$$ 
\end{cor}

\proof
By a straightforward generalization of Lemma 4.3 in \cite{ST}, which itself is
a straightforward extension of (a special case of) Lemma 4.1 from \cite{GyJS},
one can prove that if $G$ has an $s$-wide $t$-coloring and $r\ge3s-2$, then
$M_r(G)$ has an $s$-wide $(t+1)$-coloring. Thus the statement follows by
Theorem~\ref{swide} combined with the result of Stiebitz \cite{Stieb} stating
that topological $t$-chromaticity of $G$ implies topological
$(t+1)$-chromaticity of $M_r(G)$, cf. also Csorba \cite{Cs}. 
\qed

\end{document}